\documentclass[12pt,twoside,leqno]{article}
\usepackage{amsmath}
\usepackage{amssymb}
\usepackage{amsxtra}
\usepackage{amscd}
\usepackage{amsthm}
\usepackage[mathscr]{eucal}%\usepackage{eucal}

\theoremstyle{plain}
\newtheorem{thm}[subsection]{Theorem}
\newtheorem{prop}[subsection]{Proposition}
\newtheorem{cor}[subsection]{Corollary}
\newtheorem{lem}[subsection]{Lemma}

\theoremstyle{definition}

\newcommand{\n}{\noindent}
\newcommand{\hb}{\hfil\break}

\newcommand\Ker{\mathrm{Ker}\,}%
\newcommand\End{\mathrm{End}\,}%
\newcommand\Hom{\mathrm{Hom}\,}%
\newcommand\cO{\mathcal O}

\newcommand\cL{\mathcal L}

\newcommand{\Q}{{\mathbb{Q}}}
\newcommand{\Z}{{\mathbb{Z}}}
\newcommand{\R}{{\mathbb{R}}}
\newcommand{\C}{{\mathbb{C}}}

\newcommand{\upc}[1]{\overset {\lower 0.3ex \hbox{${\;}_{\circ}$}}{#1}}
\newcommand{\Gal}{\operatorname{{Gal}}}
\newcommand{\Sel}{\operatorname{{Sel}}}

\newcommand{\E}{\cal{E}}

\newcommand{\Aut}{{\rm {Aut}}}

\newcommand{\Fin}{F_{\infty}}

\newcommand{\ord}{{{\rm {ord}}}}

\newcommand{\efp}{E_{{\frak p}^{\infty}}}

\newcommand{\fp}{{\frak p}}
\newcommand{\fsp}{\frak p^*}{}
\input cyracc.def %sha
\font\tencyr=wncyr10
\def\sha{\text{\tencyr\cyracc{Sh}}}

\begin{document}

\title{The Tate-Shafarevich group for elliptic curves with
complex multiplication}

\author{J. Coates, Z. Liang, R. Sujatha}

\maketitle

%\hfill{{\it To John Cannon and Derek Holt}}

\section{Introduction}

Let $E$ be an elliptic curve over $\Q$ and put  $g_{E/\Q}={\rm rank~of~}
E(\Q)$. Let
$$
\sha(E/{\Q})\,=\, \Ker(H^1(\Q,E) \to \underset{v}{\oplus}
\,H^1(\Q_v,E)),
$$
where $v$ ranges over all places of $\Q$ and $\Q_v$ is the completion
of $\Q$ at $v$, denote its Tate-Shafarevich group. As usual,
$L(E/\Q,s)$
is the complex $L$-function of $E$ over $\Q$. Since $E$ is now
known to be modular,
Kolyvagin's work \cite{Ko} shows that $\sha(E/\Q)$ is finite
if $L(E/\Q,s)$ has a zero at $s=1$ of order $\leq 1$, and that $g_{E/\Q}$
is equal to the order of the zero of $L(E/\Q,s)$ at $s=1$. His
proof relies heavily on the theory of Heegner points and
the work of Gross and Zagier. However, when $L(E/\Q,s)$ has a zero
at $s=1$ of order $\geq 2$, all is shrouded in mystery. It is unknown
whether or not $L(E/\Q,s)$ has a zero at $s=1$ of order $\geq g_{E/\Q}$,
and no link
between $L(E/\Q,s)$ and $\sha(E/\Q)$ has ever been proven. In particular, the
finiteness of $\sha(E/\Q)$
is unknown for a single elliptic curve $E/\Q$ with $g_{E/\Q} \geq
2$. This state of affairs is particularly galling for number
theorists, since the conjecture of Birch and Swinnerton-Dyer
even gives an exact formula for the order of $\sha(E/\Q)$, which
predicts that in the vast majority of numerical examples $\sha(E/\Q)$
is zero when $g_{E/\Q}
\geq 2$. We also stress that in complete contrast to the situation for
finding $g_{E/\Q}$, it is impossible to calculate $\sha(E/\Q)$ by
classical descent methods, except for its $p$-primary subgroup for
small primes $p$, usually with $p \leq 5$.

\medskip

By contrast, in the $p$-adic world, it has long been known that the main
conjectures of Iwasawa theory provide a precise link between the
$\Z_p$-corank of the $p$-primary subgroup of $\sha(E/\Q)$ and the
multiplicity of the zero of certain $p$-adic $L$-functions at the
point $s=1$ in the $p$-adic plane, at least when $E$ has potential
good ordinary reduction at $p$. However, it seems that
little effort has been made so far to exploit this deep
connexion for theoretical purposes, and the only numerical applications
to date are given in the recent paper \cite{SW}, see also \cite{PR1},
\cite{BPR} for
the case of supersingular reduction at $p$ . The aim of this paper
is to make some modest first steps in this direction in the special case
of elliptic curves with complex multiplication. We begin with a
theoretical result. For each prime $p$, let $t_{E/\Q,p}$ denote the
$\Z_p$-corank of the $p$-primary subgroup of $\sha(E/\Q).$ While we
cannot prove the vanishing of $t_{E/\Q,p}$ for infinitely many $p$ in
any new cases, we can at least establish the following rather
general weak upper bound for $t_{E/\Q,p}$ for sufficiently large good
ordinary primes $p$.

\begin{thm}\label{t1}
Assume that $E/\Q$ admits complex multiplication. For each $\epsilon>0$,
there exists an explicitly computable number $c(E,\epsilon),$ depending
only on $E$ and $\epsilon$, such that
\begin{equation}\label{0}
t_{E/\Q,p} \leq (1+\epsilon)p-g_{E/\Q}
\end{equation}
for all primes $ p \geq c(E,\epsilon)$ where $E$ has good ordinary
reduction.
\end{thm}

We remark that a much stronger form of Theorem \ref{t1} is
known in the geometric analogue (i.e. the case of an elliptic curve
over a function field in one variable over a finite field), thanks
to the work of Artin and Tate \cite{AT}. Indeed, their work shows
that, in the geometric analogue, the number of copies of $\Q_p/\Z_p$
occuring in the Tate-Shafarevich group has an absolute upper bound
which is independent of $p$. We also note in passing that, after many
special
cases were established
by earlier authors, the Dokchitser brothers \cite{DD} have finally
proven that, for all elliptic curves $E$ over $\Q$ and all primes
$p$, the parity of $g_{E/\Q}+t_{E/\Q,p}$ is equal to the parity of
the order of zero at $s=1$ of the complex $L$-function of $E/\Q$; in
particular, the parity of $t_{E/\Q,p}$ does not depend on $p$.

\medskip

In the second part of the paper, we show that
the $p$-adic methods of Iwasawa theory enable one to
push numerical calculations of $t_{E/\Q,p}$ over a much larger range of
$p$ where $E$ admits good ordinary reduction than is possible by
classical methods. We consider the elliptic curves
\begin{equation}\label{e17}
y^2=x^3-17x
\end{equation}
and
\begin{equation}\label{e14}
y^2=x^3+14x.
\end{equation}
Both curves admit complex multiplication by the ring of Gaussian
integers
$\Z[i]$, and have $g_{E/\Q}=2.$ The conjecture of Birch and
Swinnerton-Dyer predicts that $\sha(E/\Q)=0$ for both curves.

\medskip

\begin{thm}\label{t2}
For the elliptic curves \eqref{e17} and \eqref{e14}, we have
$t_{E/\Q,p}=0$ for all primes $p$ with $p \equiv 1 \!\!\mod 4$ and
$p < 13500,$ excluding $p=17$ for \eqref{e17}. Moreover
$\sha(E/\Q)(p)=0$ for all such primes $p$.
\end{thm}

\noindent It is surprising that, for the curve \eqref{e17}, the
$p$-adic $L$-function we consider has no other zeroes beyond the
zero of order 2 arising from the fact that $E(\Q)$ has rank 2, for
all primes $p < 13500$ with $p \equiv 1 \! \mod \,4$ and $p$
distinct from 17 (more precisely, our computations show that, for
this curve and these primes $p$, the power series $H_{\frak p}(T)$
in ${\cal I}[[T]]$, whose existence is given by Proposition
\ref{hep}, is of the form $T^2.\,J_{\frak p}(T)$, where $J_{\frak
p}(T)$ is a unit in ${\cal I}[[T]])$. For the curve \eqref{e14},
there are additional zeroes for precisely the two primes $p = 29$
and $277$ amongst all $p \equiv 1 \!\!\,\mod 4$ with $p < 13500.$

\medskip

We are grateful to C. Wuthrich for his comments on our work, and some
independent numerical calculations based on \cite{SW}.

\section{$\fp$-adic $L$-functions and the main conjecture}

In this section, we briefly explain the theoretical aspects of the
Iwasawa theory of elliptic curves with complex multiplication, which
underlie the proof of Theorem \ref{t1}, and the computational work
described in \S 3. For a systematic
account of the Iwasawa theory for curves with complex multiplication,
see the forthcoming book \cite{CS}.

Let $K$ be an imaginary quadratic field, and write ${\cO}_K$ for the
ring of integers of $K$. We fix an embedding of $K$ in $\C$. Let $E$
be an elliptic curve defined over $K$ such that $\End_K(E) \otimes_{\Z}
\Q$ is isomorphic to $K$,
where $\End_K(E)$ denotes the ring of $K$-endomorphisms of $E$. It is
well-known that $E$ is isogenous over $K$ to a curve whose ring of
$K$-endomorphisms is isomorphic to ${\cO}_K$. As the results we shall
discuss depend only on the isogeny class of $E$,
we shall assume henceforth that
\begin{equation}\label{1}
\End_K(E)\simeq \cO_K.
\end{equation}
The existence of such an elliptic curve defined over $K$ implies, by the
classical theory of complex multiplication, that $K$ has class number
1. We choose a global minimal Weierstrass equation for $E$
\begin{equation}\label{2}
y^2+a_1xy+a_3y=x^3+a_2x^2+a_4x+a_6
\end{equation}
whose coefficients $a_i$ belong to ${\cO}_K$. Write $\psi_E$ for the
Gr\"ossencharacter of $K$ attached to $E$ by the theory of complex
multiplication. Recall that if $v$ is a finite place of $K$ such that
$E$ has good reduction at $v$, and if $k_v$ denotes the residue field
of $v$, then the theory of complex multiplication shows that there is
a unique element $\pi_v$ of $\End_K(E)$ such that the reduction of
$\pi_v$ modulo $v$ is the Frobenius endomorphism of the reduction of
$E$ modulo $v$, relative to $k_v$. The Gr\"ossencharacter $\psi_E$ is
then given by $\psi_E(v)=\pi_v$. We write $\frak f$ for the conductor
of $\psi_E$. It is well known that the prime factors of $\frak f$ are
precisely the primes of $K$ where $E$ has bad reduction.
For each integer $n \geq 1$, we define
$$
L_{\frak f}\,(\bar\psi_E^n,s)\,=\,\underset{(v,\frak f)=1}{\prod}\,
   \left(1-\frac{\bar\psi_E^n (v)}{(Nv)^s}\right)^{-1}.
$$
Further, $L\,(\bar\psi_E^n,s)$ will denote the primitive Hecke $L$-function of
$\bar\psi_E^n.$

Let $\cL$ be the period lattice of the N\'eron differential
$$
\varpi\,=\, \frac{dx}{2y+a_1x+a_3},
$$
and let
$$
\Phi(z,\cL):\, \C/\cL\, \simeq \, E(\Bbb C)
$$
be the isomorphism given by
$$
\Phi(z,\cL)=(\wp(z,\cL)-\frac{a_1^2+4a_2}{12},
\,\frac{1}{2}(\wp'(z,\cL)-a_1(\wp(z,\cL)-\frac{a_1^2 + 4a_2}{12})-a_3)),
$$
where $\wp(z,\cL)$ denotes the Weierstrass $\wp$-function attached to
$\cL$. Since $\cO_K$ has class number 1, there exists
$\Omega_{\infty}$ in $\C^{\times}$ such that
\begin{equation}\label{om}
\cL=\Omega_{\infty}\,{\cO}_K.
\end{equation}
As we shall explain below (see \eqref{tr}), it is well-known that
\begin{equation}\label{omp}
\Omega_{\infty}^{-n}\,L(\bar\psi_E^n,n) \,\in\,K
\end{equation}
for all integers $n \geq 1$. Moreover,
\begin{equation}\label{ln}
L(\bar{\psi}_E^n,n)\neq 0~{\rm for}~ n \geq 3,
\end{equation}
since the Euler product converges when $n \geq 3$ ( in fact, \eqref{ln}
also holds for $n=2$, but the proof is more complicated). Put
\begin{equation}\label{cpom}
c_p(E)\,=\,\Omega_{\infty}^{-p}\,L\,(\bar \psi_E^p,p).
\end{equation}

\medskip

\noindent If $\frak h$ is any integral ideal of $K$, we
define
\begin{equation}\label{15}
E_{\frak h}\,=\, \Ker \left(E(\bar K)\, \overset{h}{\longrightarrow}\,
E(\bar K)\right),
\end{equation}
where $h$ is any generator of $\frak h$.
Define $E_{{\frak p}^{\infty}}\,=\,
\underset{n\geq 1}{\cup}\,E_{{\frak p}^n}.$ Let $\cal M$ be any Galois
extension of $K$. For each non-archimedean place $w$ of $\cal M$, let
${\cal M}_w$ be the union of the completions at $u$ of all finite
extensions of $K$ contained in $\cal M$. We recall that the classical
${\frak p}^{\infty}$-Selmer group of $E$ over $\cal M$ is defined by
$$
\Sel_{\frak p}(E/{\cal M})\,=\, \Ker(H^1(\Gal(\bar{\cal M}/{\cal M}),
E_{{\frak p}^{\infty}}) \, \to \,\underset{w}{\prod}\,
H^1(\Gal(\bar{\cal M}_w/{\cal M}_w), E(\bar{\cal M}_w)),
$$
where $w$ runs over all non-archimedean places of $\cal M$. The Galois
group of $\cal M$ over $K$ operates on $\Sel_{\frak p}(E/{\cal M})$ in
the natural fashion. If $A$ is any ${\cal O}_K$-module, $A(\frak p)$
will denote the submodule consisting of all elements which are all
annihilated by some power of a generator of $\frak p$. Then we have
the exact sequence
\begin{equation}\label{24}
0 \to E({\cal M}) \otimes_{{\cal O}_K} (K_{\frak p}/{\cal O}_{\frak p})
\to \Sel_{\frak p}(E/{\cal M}) \to \sha (E/{\cal M})(\frak p) \to 0,
\end{equation}
where $\sha(E/{\cal M})$ denotes the Tate-Shafarevich group of
$E$ over $\cal M$. We will also need to consider the compact
$\Z_p$-module
\begin{equation}\label{xfp}
X_{\frak p}(E/{\cal M})\,=\, \Hom(\Sel_{\frak p}(E/{\cal
  M}),\Q_p/\Z_p).
\end{equation}
\noindent When ${\cal M}$ is any finite extension of $K$, classical
arguments from Galois cohomology show that $X_{\frak p}(E/{\cal M})$ is
a finitely generated ${\bf Z}_p$-module. In particular, we define
\begin{equation}\label{sptp}
s_{\frak p}\,=\,\Z_p{\rm -rank~of~} X_{\frak p}(E/K),~~t_{\frak
p}\,=\,\Z_p{\rm - corank~of~} \sha(E/K)(\frak p).
\end{equation}
It is clear from \eqref{24} that we have
\begin{equation}\label{stn}
s_{\frak p}\,=\,t_{\frak p} + n_{E/K},
\end{equation}
where $n_{E/K}\,=\,{\cal O}_K$-rank of $E(K)$. We denote the
number of roots of unity in $K$ by $w$.

\begin{thm}\label{cp}
Let $p$ be a prime number
such that (i) $(p,\frak f)=1,$ (ii) $(p,w)=1,$ and (iii) $p$ splits in
$K$, say $p{\cO}_K=\frak p \fsp.$ Let $m_{\fp}$ (resp. $m_{\fsp}$) denote $\ord_{\fp}(c_p(E))$
(resp. $\ord_{\fsp}(c_p(E)).$ Then we always have
\begin{equation}\label{ineq}
m_{\frak p}\geq s_{\frak p},\,\,\,m_{\fsp}\geq s_{\fsp}.
\end{equation}
Moreover, if either
$m_{\fp}= n_{E/K}$ or $m_{\fsp} =n_{E/K}$,
then  $\sha(E/K)(p)$ is finite.
\end{thm}

In fact, a stronger form
of the theorem holds if $E$ is defined over $\Q$. Assume therefore
that $E$ is defined over $\Q$, and write $L(E/\Q,s)$ for the
Hasse-Weil $L$-function of $E$ over $\Q$. By the theorem of
Deuring-Weil, we have
\begin{equation}\label{10}
L(E/\Q,s) \,=\,L(\psi_E,s),
\end{equation}
where the $L$-function on the right is the complex $L$-function
attached to the Gr\"ossencharacter $\psi_E$. Put
\begin{equation}\label{11}
g_{E/\Q}\,=\,\Z{\rm- rank}~ {\rm of~} E(\Q),\,\,
r_{E/\Q}\,=\,{\rm order~of~zero~at~} s=1~{\rm of}~ L(E/\Q,s).
\end{equation}
\noindent As $E$ is defined over $\Q$, it has real periods, and we define
$\Omega_{\infty}^+$ to be its smallest positive real period. Thus
\begin{equation}\label{omplus}
\Omega_{\infty}^+\,=\, \Omega_{\infty}\,\alpha(E),
\end{equation}
where $\alpha(E)$ is some non-zero element of ${\cal O}_K$. Put
\begin{equation}\label{cpom+}
c_p^+(E)\,=\,\left(\Omega^+_{\infty}\right)^{-p}\,L\,(\bar \psi_E^p,p).
\end{equation}
\noindent Let $\tilde E_p$ denote the reduction of $E$ modulo $p$.

\begin{thm}\label{cp1}
Assume that $E$ is defined over $\Q$. Then $c_p^+(E) \in \Q$. Let $p$ be a
prime number such that (i) $E$ has good reduction at $p$, (ii) $(p,w)=1$,
(iii) $p$ splits in $K$, and (iv) $(p,\alpha(E))=1.$ Assume also that
$r_{E/\Q}\equiv g_{E/\Q}\!\!\mod 2.$ If we have
\begin{equation}\label{cpeq}
\ord_p(c^+_p(E)) < g_{E/\Q} +2,
\end{equation}
then $\sha(E/K)(p)$ is finite.
Moreover,if
\begin{equation}\label{cpeq1}
\ord_p(c^+_p(E)) =  g_{E/\Q},
\end{equation}
and $\tilde E_p({\Bbb F}_p)$ has order prime to $p$ with $(p,6)=1$, then
$\sha(E/K)(p)=0.$
\end{thm}

\medskip

We shall say a prime $p$ satisfying (i), (ii), (iii), and (iv) of Theorem
\ref{cp1} is {\it exceptional for } $E$ if
\begin{equation}\label{exc}
\ord_p(c^+_p(E)) > g_{E/\Q}.
\end{equation}

\noindent For example,  for the curve \eqref{e14}, with
$g_{E/\Q}=2$, the primes $p=29,\,277$ are the only exceptional
primes congruent to 1 mod 4 for $p < 13500.$  However, for
$p=29,\,277$, our calculations show that $\ord_p(c_p^+(E))=3$, and
so $\sha(E/K)(p)$ is finite. For these two exceptional primes, C.
Wuthrich computed the Mazur-Swinnerton-Dyer $p$-adic $L$ function
for the curve $E$ defined by \eqref{e14}, and showed in this
way that we have also that $\sha(E/K)(p)=0$ for both primes. It is
surprising that there are no exceptional primes $p$ congruent to
1 mod 4 for the curve \eqref{e17} with $p < 13500.$

\bigskip

For all integers $n \geq 1$, let
${\E}^*_n(z,\cL)$ denote the Eisenstein series of $\cL$, as defined by
Eisenstein (see Weil \cite{We} or \cite{GS}). In particular, for $n \geq
3$, we have
\begin{equation}\label{eis}
{\E}^*_n(z,\cL)\,=\, \frac{(-1)^n}{(n-1)!}\,
\left(\frac{d}{dz}\right)^{n-2}
(\wp(z,\cL)).
\end{equation}
The following fundamental formula, which will be the basis of our
subsequent work, is proven in \cite{CW}.

\begin{thm}
Let $f$ be any generator of the conductor $\frak f$ of $\psi_E.$
Then, for all integers $n \geq 1$, we have
${\E}^*_n(\frac{\Omega_{\infty}}{f},\cL)
\in K(E_{\frak f}),$ and
\begin{equation}\label{tr}
w\,\Omega_{\infty}^{-n}\,L_{\frak f}(\bar\psi^n_E,n) \,=\,
f^{-n}\,{\rm
Trace}_{K(E_f)/K}\left({\E}^*_n\left(\frac{\Omega_{\infty}}{f},\cL\right)
\right).
\end{equation}
\end{thm}
\noindent Note that \eqref{omp} is an immediate consequence of this result.

\medskip

We  now fix a prime number $p$ satisfying $(p,\frak f)= (p,w)=1$
and $p\cO_K=\frak p\frak p^*$, where $\frak p,~\frak p^*$ are distinct
ideals  of $K$.  We pick  one of  these primes,  say $\frak  p$,  and an
embedding
\begin{equation}\label{20}
i_{\frak p}\,:\, \bar K \hookrightarrow \bar{\Q}_p,
\end{equation}
which induces $\frak p$ on $K$. For simplicity, we shall usually omit
$i_{\frak p}$ from subsequent formulae. As was shown in \cite{CW1}, (see
also \cite{CS}), there exists a $\fp$-adic  $L$-function which
essentially interpolates the image of the $L$-values \eqref{omp}.
We only state the precise result for the branch of this $\frak p$-adic
$L$-function which is needed for the proof of Theorem 1.1.

\medskip

Let $\hat E^{\frak p}$ be the formal group of $E$ at $\frak p$, so
that we can take $t=-x/y$ to be a parameter of $\hat E^{\frak p}.$ Let
$\hat{\Bbb G}_m$ be the formal multiplicative group, and write $u$ for
its parameter. Denote by $\cal I$ the ring of integers of the
completion of the maximal unramified extension of $\Q_p$. If $T$ is a
variable, then ${\cal I}[[T]]$ will denote, as usual, the ring of
formal power series in $T$ with coefficients in $\cal I$. As $\hat
E^{\frak p}$ is a formal group of height 1 ( in fact, it is even a
Lubin-Tate group over $\Z_p$ attached to the parameter $\psi_E(\frak
p)$), it is well-known that there is an isomorphism over $\cal I$
\begin{equation}\label{21}
\delta_{\frak p}\,:\, \hat{\Bbb G}_m \,\simeq\, \hat E^{\frak p},
\end{equation}
which is given by a formal power series $t=\delta_{\frak p}(u)$ in
${\cal I}[[u]].$ We can then define the $\frak p$-adic period
$\Omega_{\frak
  p}$ in ${\cal I}^{\times}$ by
\begin{equation}\label{22}
\Omega_{\frak p}\,=\,\frac{\delta_{\frak p}(u)}{u}\arrowvert_{u=0}.
\end{equation}

\medskip

\begin{prop}\label{hep}
Assume $\Omega_{\infty}$ and $\Omega_{\frak p}$ are fixed. Then there
exists a unique power series $H_{\frak p}(T)$ in ${\cal I}[[T]]$
such that, for all integers $n \geq 1$ with $ n \equiv 1 \!\!\mod
(p-1),$ we have
\begin{equation}\label{23}
\Omega_{\frak p}^{-n}\,H_{\frak p}\left((1+p)^n-1\right)\,=\,
\Omega_{\infty}^{-n}(n-1)!\,L\,(\bar
\psi_E^n,n)\,\left(1-\frac{\psi_E^n (\frak p)}{N\frak p}\right).
\end{equation}
\end{prop}
\noindent For a proof of the existence of this $\fp$-adic $L$-function
$H_{\fp}(T)$, see \cite{CW1} or
\cite{CS}. Note that when $n \equiv 1 \!\!\mod(p-1),$ $\frak f$ is the exact
conductor of $\bar \psi_E^n$, and so $L_{\frak f}\,(\bar\psi_E^n,s)$
coincides with the primitive $L$-function $L\,(\bar\psi_E^n,s).$

\medskip

This $\frak p$-adic $L$-function is related to descent theory on $E$
via the so-called ``one variable main conjecture'' for the Iwasawa
theory of $E$ over the unique $\Z_p$-extension of $K$ unramified
outside $\frak p$. Define
$$
F_{\infty}\,=\,K(E_{{\frak p}^{\infty}}), \,\,\, G\,=\,\Gal(\Fin/K).
$$
The action of $G$ on $\efp$ defines a homomorphism
\begin{equation}\label{25}
\chi_{\frak p}\,:\, G \,\to\, \Aut (\efp)=\Z_p^{\times}
\end{equation}
which is an isomorphism because $\hat{E}^{\frak p}$ is a Lubin-Tate
group.  Let $K_{\infty}$ be the unique $\Z_p$-extension contained in
$\Fin$ (class field theory shows that $K_{\infty}$ is the unique
$\Z_p$-extension of $K$ unramified outside $\frak p$). Put
$$
\Gamma\,=\,\Gal(K_{\infty}/K),\,\,\,\Lambda(\Gamma)\,=\,
\underset{\leftarrow}{\lim}\,\Z_p[\Gamma/U],
$$ where $U$ runs over the open subgroups of $\Gamma$. There is a
natural continuous action of $\Gamma$ on $X_{\frak p}(E/K_{\infty}),$
and this extends to an action of the Iwasawa algebra
$\Lambda(\Gamma)$. Since it is known that $X_{\frak p}(E/K_{\infty})$
is a finitely generated torsion $\Lambda(\Gamma)$-module (see
\cite{CW1}, \cite{CS}), it follows from the structure theory for such
modules that there is an exact sequence of $\Lambda(\Gamma)$-modules
$$
0 \to \bigoplus_{i=1}^r \,\Lambda(\Gamma)/f_i\,\Lambda(\Gamma)
\to X_{\frak p}(E/K_{\infty})\to D \to 0
$$ where $f_1,\cdots, f_r$ are non-zero elements of $\Lambda(\Gamma)$
and $D$ is a finite $\Lambda(\Gamma)$-module. We now pick the unique
topological generator $\gamma_{\frak p}$ of $\Gamma$ such that
$\chi_{\frak p}(\gamma_{\frak p})=1+p$, and write
$$
j\,:\, \Lambda(\Gamma) \, \to \, \Z_p[[T]]
$$
for the unique isomorphism of topological $\Z_p$-algebras with
$j(\gamma_{\frak p})=1+T.$ For simplicity, put
\begin{equation}\label{26}
B_{\frak p}(T)\,=\, j\,(\overset{r}{\underset{i=1}{\prod}}\,f_i).
\end{equation}
The power series $B_{\fp}(T)$ is uniquely determined up to
multiplication
by a unit in $\Z_p[[T]],$ and is called  a characteristic power series
for $X_{\fp}(E/K_{\infty}).$  We shall make essential use of the following deep result
(see \cite{Ru}, or \cite{CS}).

\begin{thm}\label{1mc}(One variable main conjecture)
$$
H_{\frak p}\left((1+p)(1+T)-1\right){\cal I}[[T]]\,=\,
B_{\frak p}(T){\cal I}[[T]].
$$
\end{thm}

\noindent In addition, we shall need (see \cite[Chap. 4, Cor. 16]{PR}):-

\begin{prop}\label{cork}
The two groups $\sha(E/K)(\frak p)$ and $\sha(E/K)(\frak p^*)$
have the same $\Z_p$-corank. In particular, one is finite if and only
if the other is also finite.
\end{prop}

\medskip

We can now prove Theorem \ref{cp}.
Since $\chi_{\fp}$ is an isomorphism, we have $\efp(K_{\infty}) =0.$
It follows that the restriction map from ${\cal
S}_{\frak p}(E/K)$ to ${\cal S}_{\frak p}(E/K_{\infty})$ is injective,
and by duality, we obtain a surjective $\Gamma$-homomorphism
\begin{equation}\label{27}
X_{\frak p}(E/K_{\infty})\, \to X_{\frak p}(E/K).
\end{equation}

\noindent Recall that $s_{\frak p}$ denotes the $\Z_p$-rank of $X_{\frak
p}(E/K)$. As
$\Gamma$ acts trivially on $X_{\frak p}(E/K)$, it follows from
\eqref{27} by a well-known property of characteristic ideals of
torsion $\Lambda(\Gamma)$-modules, that $T^{s_{\frak p}}$ must divide
$B_{\frak p}(T)$ in $\Z_p[[T]].$ Hence we conclude from Theorem
\ref{1mc} that
\begin{equation}\label{28}
H_{\frak p}\left((1+p)(1+T)-1\right)\,=\,T^{s_{\frak p}}\,h(T)
\end{equation}
for some $h(T)$ in ${\cal I}[[T]]$. Evaluating both sides at
$(1+p)^{n-1}-1$ for any $n$ in $\Z$, we conclude that we always have
\begin{equation}\label{29}
H_{\frak p}\left((1+p)^n-1\right)\,\equiv \,0\!\!\mod p^{s_{\frak p}}.
\end{equation}
Taking $n=p,$ and noting that $\left(1-\frac{\psi_E(\frak
  p)^p}{N{\frak p}}\right)$ is a
unit at $\frak p$, we conclude from \eqref{29} and Proposition
\ref{hep} that
\begin{equation}\label{30}
c_p(E) \equiv 0\!\!\mod {\frak p}^{s_p}.
\end{equation}
Replacing $\frak p$ by $\frak p^*,$ the same argument shows that
\begin{equation}\label{32}
c_p(E) \equiv 0\!\!\mod (\frak p^{*})^{s_{\frak p^{*}}}.
\end{equation}
Hence \eqref{ineq} follows. Moreover, if $m_{\fp}=n_{E/K}$, then
$t_{\fp}=0$ and so $t_{\fsp}=0$ by Proposition \ref{cork}.
A similar argument holds if $m_{\fsp}=n_{E/K}$. This completes the
proof of Theorem \ref{cp}.
\qed

\medskip

\begin{cor}
We have $m_{\fp}=n_{E/K}$ if and only if the characteristic power series
of $X_{\fp}(E/K_{\infty})$ can be taken to be $T^{n_{E/K}}.$
\end{cor}
\proof
If $m_{\fp}=n_{E/K}$, the above argument shows that we must have
$s_{\fp}=n_{E/K},$ and $h(0)$ a $p$-adic unit. It follows from Theorem
\ref{1mc} that $B_{\fp}(T)$ must be of the form $T^{n_{E/K}}$ times a
unit in $\Z_p[[T]]$. Conversely, if the characteristic power series of
$X_{\fp}(E/K_{\infty})$ can be taken to be $T^{n_{E/K}},$ then Theorem
\ref{1mc} shows that $H_{\fp}(T)$ is equal to $T^{n_{E/K}}$ times a unit
in ${\cal I}[[T]],$ whence it is plain that $m_{\fp}=n_{E/K}$. This
completes the proof.
\qed

\medskip

 Our numerical calculations show that, for the elliptic curve
$$
y^2=x^3-17x, ~{\rm with}~ n_{E/K}=2,
$$
we have $m_{\fp}=2$ for all primes $p$ with $p\equiv 1 \!\!\mod 4$,
$p\neq 17,$ and $p<13500.$ Thus the characteristic power series of
$X_{\fp}(E/K_{\infty})$ is $T^2$ for all such primes. On the other
hand, for the elliptic curve
$$
y^2=x^3+14x,~{\rm with}~ n_{E/K}=2,
$$
we have $m_{\fp}=2$ for all primes $p$ with $p \equiv 1\!\!\mod 4$
and $p < 13500$, except $p=29,\,277.$ Thus for all such primes, with
the exception of these two, the characteristic power series of
$X_{\fp}(E/K_{\infty})$ is $T^2$.

\medskip

We now establish Theorem \ref{cp1}. Assuming that $E$ is defined over
$\Q$, we have
\begin{equation}\label{psif}
\bar{\frak f} \,=\,\frak f,~{\rm and}~\psi_E(\bar{\frak
a})\,=\,\overline{\psi_E(\frak a)}
\end{equation}
for all integral ideals $\frak a$ of $K$ with $(\frak a,\frak f)=1.$
Hence
$$
L(\psi_E^p,s)\,=\,L(\bar \psi_E^p,s).
$$
Evaluating at $s=p,$ we conclude that
$$
L(\bar \psi_E^p,p)\in \R.
$$
As $\Omega^+_{\infty}$ is real, it follows that
$$
c^+_p(E) \in K \cap \R = \Q.
$$

As before, let $t_{\fp}$ (resp. $t_{\fsp})$ be the $\Z_p$-corank of
$\sha(E/K)(\fp)$ (resp. $\sha(E/K)(\fsp)$), and let $t_{E/\Q,p}$ be the
$\Z_p$-corank of $\sha(E/\Q)(p).$ Then we claim that
\begin{equation}\label{utp}
t_{E/\Q,p}=t_{\fp}=t_{\fsp}.
\end{equation}
Indeed, the second equality is just Proposition \ref{cork}. To prove
the first equality, note that $E$ is isogenous over $\Q$ to the twist
$E'$ of $E$ by the quadratic character of $K$ (see, for example,
\cite{Gr}). Thus, $\sha(E/\Q)(p)$ and $\sha(E'/\Q)(p)$ have the same
$\Z_p$-corank, and hence the $\Z_p$-corank of $\sha(E/K)(p)$ is equal
to $2 t_{E/\Q,p}$. On the other hand, the $\Z_p$-corank of $\sha(E/K)(p)$
is
clearly equal to $t_{\fp}+t_{\fsp}=2t_{E/\Q,p}$, by Proposition
\ref{cork}. Hence $t_{E/\Q,p}=t_{\fp}$, thereby proving \eqref{utp}.

Assume now that $\sha(E/K)(p)$ is infinite, so that $t_{E/\Q,p}>0$. The
parity theorem for $E/\Q$ and the prime $p$ (due to Greenberg in
this case, but see the more
general results of \cite{DD}, \cite{Ne}) asserts that
$$
g_{E/\Q}+t_{E/\Q,p}\equiv r_{E/\Q}\!\!\mod 2.
$$
By our hypothesis that $g_{E/\Q}$ and $r_{E/\Q}$ have the same
parity, it follows that $t_{E/\Q,p}$ must be even, and therefore
$t_{E/\Q,p} \geq 2$, in particular. Hence by \eqref{utp} $t_{\fp} \geq 2$.
Noting that $g_{E/\Q}=n_{E/K}$, and that $(p,\alpha(E))=1$, we conclude
from \eqref{30} that
$$
\ord_p(c^+_p(E)) \geq g_{E/\Q} +2.
$$
Hence, if \eqref{cpeq} holds, then we must have $\sha(E/K)(p)$ is
finite.

Assume now that $\ord_p(c^+_p(E))=g_{E/\Q}$. We
deduce easily from Theorem \ref{1mc} and \eqref{28}, that
$$
B_{\fp}(T)\,=\, T^{n_{E/K}}R_{\fp}(T),
$$
where $R_{\fp}(T)$ is a unit in $\Z_p[[T]]$, so
that $R_{\fp}(0)$ is a unit in $\Z_p.$
Hence, by an important
general theorem of Perrin-Riou \cite{PR}, it follows that
the canonical $\fp$-adic height pairing
$$
<~,~>_{\fp}\,:\, E(K)\otimes_{\cO} \Z_p \, \times\, E(K)
\otimes_{\cO} \Z_p \,\to \,\Q_p,
$$
where $\cO$ is embedded in $\Z_p$ via $i_{\fp}$, is non-degenerate.
Further, we have that
\begin{equation}\label{htp}
\#(\sha(E/K)(\fp)) \times \det <~ ,~>_{\fp} \times
\left(1-\frac{\psi_{E/K}(\fp)}{N{\fp}}\right)
\end{equation}
is also a $p$-adic unit, where $\det$ denotes the determinant of the
height pairing;
for this last assertion, we need our
hypothesis that $(p,6)=1.$  However, if $\tilde E_p({\Bbb F}_p)$ has
order prime to $p$ and $(p,6)=1,$ then it follows from the results of
\cite{PR} that $\det\,<~,~>_{\fp}$ is a $p$-adic integer. Hence we
conclude from \eqref{htp} that $\sha(E/K)(\fp)$ is trivial. A similar
argument proves the corresponding statement for $\sha(E/K)(\fsp)$ and
this completes the proof of Theorem \ref{cp1}.
\qed

\medskip

We next establish an upper bound for $t_{\fp}$ and $t_{\fsp}$ when $p$ is
a sufficiently large prime which splits in $K$ as $p{\cal O}_K=\fp\fsp.$

\begin{thm}\label{tbd}
For each $\epsilon >0$, there exists an explicitly computable number
$c(E,\epsilon),$ depending only on $E$ and $\epsilon$, such that
\begin{equation}\label{bdsp}
t_{\frak p}  \leq (1+\epsilon)p-n_{E/K},~~~~t_{\fsp}\leq
(1+\epsilon)p-n_{E/K},
\end{equation}
for all primes $p \geq c(E,\epsilon)$ which split in $K$ as $p{\cal
O}_K=\fp\fsp.$
\end{thm}

\noindent We note that, when $E$ is defined over $\Q$, Theorem \ref{t1} is
an immediate consequence of this result, since, thanks to \eqref{utp}, we
then have $t_{E/\Q,p}=t_{\fp},~~n_{E/K}=g_{E/\Q}.$

\medskip

We now give the proof of Theorem \ref{tbd} which is a simple
application of the formula \eqref{tr}, and the fact that $L(\bar
\psi_E^p,p)\neq 0$ (recall that the latter assertion is true because
the Euler product for $L(\bar \psi_E^p,s)$ converges for $s=p$). Put
\begin{equation}\label{thp}
\Theta_p\,=\, {\rm Trace}_{K(E_{\frak
f})/K}\left({\E}_p^*\left(\frac{\Omega_{\infty}}{f},\cL \right)\right).
\end{equation}
\noindent We emphasize that in the proof $E$ is fixed and $p$ is varying
over all sufficiently large prime numbers which split in $K$.

\begin{lem}\label{thb}
We have $\left|\Theta _p \right| \leq d_1^p,$ where $d_1>1$ is a real
number depending only on $E$ and not on $p$.
\end{lem}
\proof
We may assume $p \geq 3$.
By \eqref{eis}, we have
\begin{equation}\label{eis1}
{\E}^*_p\left(\frac{\Omega_{\infty}}{f},{\cL}\right) \,=\,
\frac{(-1)^p}{(p-1)!}
\,\left(\frac{d}{dz}\right)^{p-2}\,\left(\wp(z,\cL)\right)|_{
z=\frac{\Omega_{\infty}}{f}}.
\end{equation}
Let $\cal B$ denote a set of integral ideals of $K$, prime to $\frak f$,
such that the Galois group of $K(E_{\frak f})/K$ consists precisely of the
Artin symbols $\sigma_{\frak b}$ of the ideals $\frak b$ in $\cal B$. From
the definition of the Gr\"ossencharacter $\psi_E$ and \eqref{eis1}, we
have
$$
{\E}^*_p\left(\frac{\Omega_{\infty}}{f},{\cL}\right)^{\sigma_{\frak b}} \,
= \,{\E}^*_p\left(\psi_E(\frak b)\,\frac{\Omega_{\infty}}{f},{\cL}\right).
$$
Thus, by Cauchy's integral formula, we obtain
$$
{\E}^*_p\left(\frac{\Omega_{\infty}}{f},{\cL}\right)^{\sigma_{\frak b}} \,
=\, \frac{(-1)^p}{(p-1)\,2\pi i}
\int_{{\cal C}_{\frak b}}\,\frac{\wp(z,{\cal
L})dz}{\left(z-\frac{\psi_E(\frak b)\Omega_{\infty}}{f}\right)^{p-1}},
$$
where ${\cal C}_{\frak b}$ is a circle with centre $\frac{\psi_E(\frak
b)\Omega_{\infty}}{f}$ and sufficiently small radius so that no element of
${\cL}$ lies in or on ${\cal C}_{\frak b}$. Estimating the integral, it is
plain that
$$
\left |{\E}^*_p\left(\frac{\Omega_{\infty}}{f},{\cL}\right)^{\sigma_{\frak
b}}\right | \leq d_2^p,
$$
where $d_2>1$ depends only on $E$. Summing over all $\frak b$ in ${\cal
B}$, the assertion of the lemma follows.
\qed

\begin{lem}\label{d3}
There exists a rational integer $d_3 >1, $ depending only on $E$
and not on $p$, such
that
$$
d_3^p\,(p-1)! \,{\E}_p^*\left(\frac{\Omega_{\infty}}{f},\cal L\right)
$$
is an algebraic integer.
\end{lem}
\proof
We may assume that $p \geq 5$. Since
$$
{\E}_p^* (\lambda z,\lambda {\cL})\,=\,\lambda^{-p}\,{\E}_p^* (z,{\cL})
$$
for any complex number $\lambda,$ it suffices to prove the lemma
when our generalized
Weierstrass equation \eqref{2} for $E$ has the property that
$g_2(\cL)/2$ and $g_3(\cL)$
both belong to ${\cO}_K$; here $g_2(\cL)$ and $g_3(\cL)$ denote
the usual Weierstrass
invariants attached to \eqref{2}. Now the differential equation
$$
\left( \wp'(z,\cL)\right)^2\,=\,4\wp(z,\cL)^3-g_2(\cL)\wp(z,\cL) -g_3(\cL)
$$
implies that
$$
\wp^{(2)}(z,\cL)\,=\,6\wp(z,\cL)^2-\frac{g_2(\cL)}{2}.
$$
A simple recurrence argument on $n$ then shows that,
for all $n \geq 1$, we have
$$
\wp^{(2n)}(z,\cL)\,=\,D_n(\wp(z,\cL)),
$$
where $D_n(X)$ is a polynomial in ${\cO}_K[X]$ of degree $n+1$.
It follows immediately that
$$
\wp^{(2n+1)}(z,\cL)\,=\,B_n(\wp(z,\cL))\wp'(z,\cL),
$$
where $B_n(X)=\frac{d}{dX}(D_n(X))$ is a polynomial of degree $n$
in $\cO_K[X].$ Taking $n=\frac{p-3}{2}$, the
assertion the lemma is now clear from \eqref{eis1}, on taking $d_3$ to be
a positive integer such that
$$
d_3.\wp(\frac{\Omega_{\infty}}{f},\cL),\,\,d_3.\wp'
(\frac{\Omega_{\infty}}{f},\cL)
$$
are algebraic integers.
\qed

\medskip

We can now complete the proof of Theorem \ref{tbd}. We may assume
that $(p,\frak f)=(p,w)=1.$ By
\eqref{tr}, we then have
\begin{equation}\label{oml}
\left|\Omega_{\infty}^{-p}\,L(\bar{\psi}_E^p,p)\right|_{\frak p}\,=\,
\left|\Theta_p\right|_{\frak p}\,=\,\left|(p-1)!\,
\Theta_p\right|_{\fp},
\end{equation}
and similarly for $\fsp$.
Moreover, in view of Lemmas \ref{thb} and \ref{d3},
$$
d_3^p(p-1)!\,\Theta_p
$$
is an element of ${\cO}_K$ whose complex absolute value is at most
$d_4^p\,(p-1)!,$ where $d_4>1$ does not depend on $p$.
Since $\Theta_p\neq 0$ because $L(\bar\psi^p_E,p)\neq 0$,
we conclude from the product formula that
\begin{equation}\label{d4}
\left|d_3^p\,(p-1)!\Theta_p\right|_{\frak p}\times
\left|d_3^p(p-1)!\Theta_p\right|_{\fsp} \geq
d_4^{-2p}\left((p-1)!\right)^{-2}.
\end{equation}
It follows that, we conclude that for each $\epsilon>0$, we have
$$
\left|\Theta_p\right|_{\fp} \times \left|\Theta_p\right|_{\fsp}
\geq p^{-2(1+\epsilon)p}
$$
for all $p\geq c(E,\epsilon)).$ On the other hand, by Theorem
\ref{cp},
and \eqref{oml},
we have
$$
\left|\Theta_p\right|_{\fp} \times \left|\Theta_p\right|_{\fsp}
\leq p^{-(s_{\fp} +s_{\fsp})}.
$$
Thus
$$
s_{\fp} +s_{\fsp} \leq 2(1+\epsilon)p
$$
when $p \geq c(E,\epsilon).$ As $s_{\fp}=s_{\fsp}$, the proof of
the theorem is complete.
\qed

\medskip

Define the $p$-adic $L$-functions
$$
{\frak L}_{E,\fp}(s)\,=\, H_{\fp}\left((1+p)^s-1\right),
\,\,\,{\frak L}_{E,\fsp}(s)\,=\,H_{\fsp}\left((1+p)^s-1\right),
$$
where $s$ is now a variable in $\Z_p$. Put
$$
r_{E,\fp}\,=\,\ord_{s=1}{\frak L}_{E,\fp}(s),\,\,\,
r_{E,\fsp}\,=\,\ord_{s=1}{\frak L}_{E,\fsp}(s).
$$
We end this section by remarking that exactly the same argument which
establishes Theorem \eqref{tbd} proves the following result.
\begin{thm}\label{refp}
For each $\epsilon >0$, there exists an explicitly computable number
$c(E,\epsilon)$ such that
$$
r_{E,\fp}  +  r_{E,\fsp} \leq 2 (1+\epsilon)p
$$
for all primes $p \geq c(E,\epsilon)$ with $p$ splitting in $K$ as
$p\cO_K=\fp\fsp.$
\end{thm}

\section{Computations for $y^2=x^3-Dx$}

The goal of this section is to explain how one can use formula
\eqref{tr} to compute
$$
c_p(E)\,=\,\Omega_{\infty}^{-p}\,L(\bar \psi_E^p,p)
$$
in practice, for the family of curves
$$
E\,:\,y^2\,=\, x^3-Dx,
$$
where $D$ is a fourth-power free non-zero rational integer. For this
family of curves, $K=\Q(i)$
and the isomorphism \eqref{1} is given explicitly by mapping $i$ to
the endomorphism which sends $(x,y)$ to $(-x,iy).$  See \cite{BGS}, 
\cite{FK} for earlier computational work on the Iwasawa theory of this
family of curves.

We begin by analysing the Galois theory of the fields $K(E_{\frak f})$
where $\frak f$ again denotes the conductor of $\psi_E$. If $\frak h$
is any integral ideal of $K$, we write
$$
\phi(\frak h) \,=\,\#\left(\left(\Z[i]/{\frak h} \right)
^{\times}\right).
$$
The next lemma is a very easy
consequence of  the existence of the  Gr\"ossencharacter $\psi_E$ (see
\cite{CW},\,Lemma 3, or \cite{C},\,Lemma 7) and the fact that no root
of unity in $K$ is $\equiv 1 \!\!\mod {\frak h}$, when $\frak h$ is a
multiple of $\frak f$.

\begin{lem}\label{rcf}
Let $\frak h$  be any integral ideal of $K$ which  is divisible by the
conductor $\frak f$ of  $\psi_E$. Then $K(E_{\frak h})$ coincides with
the ray class field of $K$ modulo $\frak h$.
In particular, the degree of $K(E_{\frak h})/K$ is equal to
$\phi(\frak h)/4.$
\end{lem}

\medskip

The following well-known lemma computes $\frak f$ for the curve $E$.

\begin{lem}\label{delta}
Let $\Delta$ be the product of the distinct prime factors of $D$. Then
$\frak f=4\Delta\,\Z[i]$ if $D \not\equiv 1\!\!\mod 4$ and $\frak f=
(1+i)^3\Delta\,\Z[i]$ if $D \equiv 1 \!\!\mod 4.$
\end{lem}

\medskip

Let $E'$ denote the elliptic curve in our family with $D=1$,  i.e.
\begin{equation}\label{e'}
E'\,:\, y^2\,=\,x^3-x.
\end{equation}

\begin{lem}\label{ee'}
Assume that $E=E^D$ with $D$ divisible by an odd power of an odd
  prime. Then the extension
$K\left(E_{(1+i)^k}\right)$ is equal to $K$
  when $k=1$, to $K(D^{1/2})$ when $k=2$, and to $K(D^{1/4})$ when
  $k=3$. For $k \geq 3$, we have
$$
K\left(E_{(1+i)^k}\right) \,=\,K\left(D^{1/4},E'_{(1+i)^k}\right)
$$
and this field has degree $2^{k-1}$ over $K$, and degree 4 over
$K\left(E'_{(1+i)^k}\right).$
\end{lem}
\proof
The assertions for $k=1$ and $k=2$ are readily verified. Put
$\alpha=D^{1/4}.$ Over $K(\alpha)$, we have an isomorphism
\begin{equation}\label{iso}
E \simeq E'
\end{equation}
given by mapping the point $(x,y)$ on $E$ to the point $\left(
x/{\alpha}^2,y/{\alpha}^3\right)$ on $E'.$ Now $E'$ has conductor
$(1+i)^3,$ and $K\left( E'_{(1+i)^3} \right) =K,$ whence it follows
from \eqref{iso} that $K\left( E_{(1+i)^3} \right)=K(\alpha).$
Similarly, if  $k\geq 3$, then \eqref{iso} implies that $K
\left( E_{(1+i)^k} \right)=K\left(\alpha,E'_{(1+i)^k}\right)$. Now
Lemma \ref{rcf} applied to $E'$ shows that the degree of
$K\left(E'_{(1+i)^k}\right)$ over $K$ is $2^{k-3}$ when $k \geq
3$. Moreover, as $E'$ has good reduction outside the prime
$(1+i){\Z[i]}$, this is the only prime of $K$ which can ramify in the
extension $K\left(E'_{(1+i)^k}\right).$ Hence  $[K(\alpha):K]=4$,
and $K(\alpha)\cap K\left(E'_{(1+i)^k}\right)=K$
because of the existence of the odd prime factor dividing $D$ to an odd
power.
This completes the proof of the lemma.
\qed

%\begin{lem}\label{lD}
%Assume that $D$ is odd. Then the degree of $K(E_D)/K$ is either
%$\phi(D\Z[i])$ or $\phi(D\Z[i])/2.$
%\end{lem}
%\proof
%By the Weil pairing, $K(E_D)$ contains the field generated over $K$ by
%the
%$|D|$-th roots of unity. Hence $K(E_D)$ contains $\sqrt{D}$
%(the sign of $D$ is irrelevant since $K$ contains the fourth roots of
%unity). As above,
%let $\alpha=D^{1/4}$. Then
%$K(E_D,\alpha)=K(E'_D,\alpha)$, where $E'$ is the curve \eqref{e'},
%because $E$ is isomorphic to $E'$ over $K(\alpha)$. But the degree
%$[K(E'_D):K]$ is equal to $\phi(D{\Bbb Z}[i]),$ because $D$ is odd
%(for example $K(E'_D)=K(E'_{D(1+i)^3}),$ and apply Lemma \ref{rcf} to
%$E').$ Also $\alpha^2$ belongs to $K(E_D),$ since this latter field
%also
%contains the $|D|$-th roots of unity. Thus $K(E_D,\alpha)$ has degree
%$r\phi(D\Z[i])$ over $K$, where $r=1$ or $2$, proving the lemma.
%\qed

\begin{lem}\label{galois}
Assume that $D$ is odd. Then the degree of $K(E_D)/K$ is $\phi(D
\Z[i]).$
\end{lem}
\proof
We can assume $D\neq 1$.
By the Weil pairing, $K(E'_D)$ contains the field generated over $K$ by
the $|D|$-th roots of unity. Hence $K(E'_D)$ contains $\sqrt{D}$
(the sign of $D$ is irrelevant since $K$ contains the fourth roots of
unity).  As above, let $\alpha=D^{1/4}$. Thus $K(E'_D,\alpha)$ has
degree at most  2 over $K(E'_D)$.

Let $R_D$ denote the ray class field of $K$
modulo $D\Z[i]$. Let $(u,v)$ be a  primitive $D$-division point on
$E$. Then the classical theory of complex multiplication shows that
$R_D=K(u^2),$ and that $[R_D:K]=\phi(D\Z[i])/4.$ To prove the lemma,
it therefore suffices to show that there exists an element
$\tau$ of $\Gal(K(E_D)/K)$ such that $\tau$ fixes $R_D$, and $\tau$ is
of exact order 4. We do this as follows. As remarked in the previous
paragraph, $K(E_D')$
has degree $\phi(D\Z[i])$ over $K$ because $D$ is odd. Moreover, a
primitive $D$-division point on $E'$ is given by $(u',v')$, where
$u'=u/\alpha^2,$ $v'=v/\alpha^3.$ Recalling that multiplication by $i$
on $E'$ is given by sending $(x,y)$ to $(-x,iy)$, it follows that
there exists $\sigma$ in $\Gal(K(E'_D)/K)$ such that
\begin{equation}\label{sig}
\sigma(u',v')\,=\,(-u',iv').
\end{equation}
Now let $\sigma$ denote any extension of $\sigma$ to the field
$K(E'_D,\alpha) \,=\,K(E_D,\alpha).$ Since this field has degree at
most 2 over $K(E'_D),$ we must have that either
$\sigma(\alpha)\,=\,-\alpha$ or $\sigma(\alpha)=\alpha.$
Applying $\sigma$ to $(u',v'),$ we conclude from \eqref{sig} that
$$
\sigma\,u\,=\,u,\,\sigma\,v\,=\,v~{\rm
  or}~\sigma\,u\,=\,-u,\,\sigma\,v\,
=\,iv.
$$
It follows
from these formulae that $\sigma^4$ fixes $K(E_D),$ but $\sigma^2$
does not. Also $\sigma$ fixes $R_D$. Hence we may take $\tau$ to be
restriction of $\sigma$ to $\Gal(K(E_D)/K)$, and the proof of the
lemma is complete.
\qed

\begin{lem}\label{2D}
Let $D=2^aM$, where $a=1$ or $3$, and $M$ is odd. Then $K(E_M)$ has degree $\phi(M\Z[i])$ over
$K$, and $K(E_M,D^{1/4})$ has degree $4 \phi(M\Z[i])$ over $K$.
\end{lem}
\proof
As remarked earlier, $K(E'_M)$  has degree
$\phi(M \Z[i])$ over $K$ because $M$ is odd. Also, by Lemma
\ref{rcf}, $K(E'_{8M})$ is the ray class field of $K$ modulo $8M$,
and hence we have
$$
[K(E'_{8M}):K]\,=\,8 \phi(M\Z[i]).
$$
Since $[K(E'_8):K]=8$ by Lemma \ref{rcf}, we conclude that
\begin{equation}\label{kem}
K(E'_M) \cap K(E'_8)\,=\, K.
\end{equation}
By the Weil pairing, $K(E'_M)$ contains the field of
of $|M|$-th roots of unity, and hence also $\sqrt{M}.$
Similarly, $K(E'_8)$ contains the eighth roots of unity, and so
also $\sqrt{2}$. But $K(\sqrt{2})/K$ is an extension of degree 2,
and thus, by \eqref{kem}, $\sqrt{2}$ does not belong to $K(E'_M)$.
It follows that $\sqrt{D}$ does not belong to $K(E'_M)$ since
$a=1$ or $3$. Hence
$$
[K(E'_M,\alpha):K]\,=\, 4 \phi(M \Z[i]).
$$
But $E$ and $E'$ are isomorphic over $K(\alpha),$
whence
$$
K(E'_M,\alpha)\,=\,K(E_M,\alpha).
$$
On the other hand, it is clear that $[K(E_M,\alpha):K]$
divides $4\phi(M \Z[i]).$ It follows that
$$
[K(E_M):K]=\phi(M\Z[i]),\,\,\,[K(E_M,\alpha):K(E_M)]=4,
$$
and the proof of the lemma is complete.
\qed

\medskip

We now briefly describe the theoretical steps underlying our numerical
calculations of $\ord_p(c_p^+(E))$ for the curve $E$ when $D$ is
divisible by at least one odd prime. The Weierstrass equation
associated to $E$ is
\begin{equation}\label{weier}
\wp'(z,\cL)^2\,=\, 4 \wp(z,\cL)^3 - 4D\wp(z,\cL).
\end{equation}
Write $\frak f= f\Z[i]$ for the conductor of $\psi_E$, and define
\begin{equation}\label{up}
u\,=\,\wp\left(\frac{\Omega_{\infty}}{f}\,,\, \cL\right),\,\,\,
v= \left(\wp'\left(\frac{\Omega_{\infty}}{f}\,,\, \cL\right)\right)/2.
\end{equation}
By Lemma \ref{rcf}, $K(E_{\frak f})$ is the ray class field of $K$
modulo $\frak f$. Hence
\begin{equation}\label{kef}
K(E_{\frak f})\,=\, K(u^2)\,=\,K(u),\,\,\,v\in K(u),
\end{equation}
and the degree of $K(E_{\frak f})$ over $K$ is $d=\phi(\frak f)/4.$
As $\frak f$ is divisible by at least two distinct primes of $K$, a
theorem of Cassels \cite{Ca} shows that both $u$ and $v$ are algebraic
integers. Moreover, we can compute explicitly the monic irreducible
polynomial of $u$ over $\Z[i]$, which has degree $d$, and which we
denote by $G(X)$. Once we have computed this polynomial $G(X)$, we can
determine
\begin{equation}\label{sm}
s_m\,=\, {\rm Trace}_{K(E_{\frak f})/K}\,(u^m)~~~~(m=1,2,\cdots,d-1)
\end{equation}
recursively, using the following classical formula. Let
$$
G(X)\,=\,(X-u_1)....(X-u_d)\,=\,
X^d-\sigma_1X^{d-1}+ \cdots + (-1)^d\sigma_d
$$
where $\sigma_1,\cdots , \sigma_d$ are the elementary symmetric
fuctions in $u_1,\cdots ,u_d.$ Then we have (see for example,
\cite[Vol I, p.81]{vdw}),
\begin{equation}\label{sm1}
s_m\,=\,(-1)^{m-1}m\sigma_m + \overset{m-1}{\underset{h=1}{\Sigma}}\,
(-1)^{h-1}\,s_{m-h}\sigma_h\,\,\,(m \leq d)
\end{equation}

Now we recall that, by virtue of formulae \eqref{eis} and \eqref{tr},
we have
\begin{equation}\label{cp2}
c_p^+(E)\,=\,
-w^{-1}(f\alpha(E))^{-p}\left((p-1)!\right)^{-1} \Xi_p,
\end{equation}
where $\alpha(E)$ is as in \eqref{omplus}, and
$$
\Xi_p\,=
\,{\rm
  Trace_{K(E_{\frak f})/K}}
\left(\wp^{(p-2)}\left(\frac{\Omega_{\infty}}{f}, \cL\right)\right)
$$
for all odd primes $p$. For our curve $E=E_D$, we have $w=4$. Moreover,
we have
$$
\Omega_{\infty}^+\,=\, \Omega/\,D^{1/4}~~{\rm if}~ D>0,~~
\Omega_{\infty}^+\,=\, \Omega/\,(-D/4)^{1/4}~~{\rm if}~ D<0,
$$
where $\Omega=2.622058\cdots$ is the least positive real period
of the curve $E'$ \eqref{e'}.
Hence $\alpha(E)\,=\,1$ when $D>0$, and
$\alpha(E)\,=\,(1+i)$ when $D<0$.
 We now fix the value of $f$ following the four cases:- (i) $D>0$ and $D
\equiv 1 \!\mod\,4$, (ii) $D<0$ and $D \equiv 1 \!\mod\,4$, (iii) $D>0$ and $D
\not\equiv \! 1 \mod \,4$, and (iv) $D<0$, and $D \not\equiv \! 1 \mod \,4$.
Following these four cases, we take $f$ to be $2(1+i)\Delta,$
$(1+i)^3\Delta,$ $4\Delta$, and $4\Delta,$ so that the respective
values of $f\alpha(E)$ are given
by $2(1+i)\Delta,$ $-4\Delta$, $4\Delta$, and $4\Delta(1+i)$.

As explained in the proof of Lemma
\ref{d3}, we have
\begin{equation}\label{pp2}
\wp^{(p-2)}\left(\frac{\Omega_{\infty}}{f}\,,\,\cL \right)\,=\, B_{\frac{p-3}{2}}
\left(\wp\left(\frac{\Omega_{\infty}}{f}\,,\,\cL \right)\right)
\wp'\left(\frac{\Omega_{\infty}}{f}\,,\,\cL \right),
\end{equation}
where $B_{\frac{p-3}{2}}(X)$ is a polynomial in $\Z[X]$ of degree
$(p-3)/2$. This polynomial can easily computed recursively, using
the differential equation \eqref{weier} (see the explicit examples
below when $D=17$ and $D=-14$).

As the theory tells us that $v \in K(u),$ there exists a
polynomial $J(X)$ in $K[X]$ such that
\begin{equation}\label{pinf}
\wp'\left(\frac{\Omega_{\infty}}{f}\,,\,\cL \right)
\,=\, J\left(\wp\left(\frac{\Omega_{\infty}}{f}\,,\,\cL \right)\right).
\end{equation}

\medskip

\noindent In fact, in the numerical examples we have
considered, it is always the
case that $J(X)$ belongs to $\Z[i][1/f][X]$, and we shall assume
henceforth that this is the case. Hence, multiplying together
$B_{\frac{p-3}{2}}(X)$ and $J(X)$, and using
the fact that $G\left(\wp\left(
\frac{\Omega_{\infty}}{f}, \cal L\right)\right)=0$, we deduce that
$$
\wp^{(p-2)}\left(\frac{\Omega_{\infty}}{f}, \cal L\right)\, =
\, A_p\left(\wp\left(\frac{\Omega_{\infty}}{f}, \cal L\right)\right),
$$
where $A_p(X)$ is a polynomial in $\Z[i][1/f][X]$ of degree at most
$d-1$. Writing $A_p(X)\,=\,\overset{d-1}{\underset{j=0}{\Sigma}}\,
a_{j,p} X^j,$ it follows that
$$
\Xi_p\,=\, \overset{d-1}{\underset{j=0}{\Sigma}}\,
a_{j,p} s_j,
$$
and we can then compute $c_p^+(E)$ using the formula \eqref{cp2}. The
machine then calculates $\ord_p(c_p^+(E))$ (which our theory shows is
always
$\geq 0$), followed by
$$
c_p^+(E)\!\mod p^k,~ {\rm where}~ k=\ord_p(c_p^+(E))+1.
$$
\smallskip
\noindent Finally, we note that $\tilde E_p({\Bbb F}_p)$ has order prime
to $p$ for all $p > 5$
with $(p,D)=1$. This is clear when $p \equiv 3 \!\mod 4$, since
then $\tilde E_p$ is supersingular. For $p\equiv 1 \! \mod 4$, say
$p{\Bbb Z}[i] = \frak p. \frak p*$, we have
$a_p = {\rm Trace}
_{K/{\Bbb Q}}(\psi_E(\frak p))$ must be even
because
2 is ramified in $\Q(i)$, from which it follows that we cannot have
$a_p = 1$, which clearly implies the assertion for these primes $p$.

\medskip

The computations described above have been carried out for the two
curves $D=17$ and $D=-14$ for all primes $p$ with $p \equiv \!1 \mod 4$
and $p < 13,500$ (the prime $p=17$ is excluded when $D=17$). We have
$$
\begin{array}{lcl}
D=17,& f=2(1+i)17,& d=256\\
D=-14, & f=56, & d=384.
\end{array}
$$
\noindent For both cases, the polynomials $ G(X),\,J(X),\,
B_{\frac{p-3}{2}}(X),\, A_p(X)$ have been computed explicitly, and
are given at \cite{tab}, as they are too elaborate to include here.
However, as an illustrative example where the coefficients are still
not too enormous, we give now the polynomials $B_{13}(X)$, which
occur for $p=29,$
$$
\begin{array}{lc}
D=-14 &\\
B_{13}(X)= 7496723869173\times 2^{24} (431525237696 X + 3877463640960
X^3 +
   5545863414000 X^5 & \\
   + 2565173520000 X^7 + 490959787500 X^9 +
   40724775000 X^{11} + 1212046875 X^{13})& \\
   D=17&\\
 B_{13}(X)=
   7496723869173\times 2^{24} (1383348216959 X - 10236515835780 X^3 +
   12057373443375 X^5 & \\
    - 4592819790000 X^7 + 723915196875 X^9 -
   49451512500 X^{11} + 1212046875 X^{13})& .
   \end{array}
$$
For these two curves, and our range of $p$, our computations show that
$\ord_p(c_p^+(E))=2,$ except for the two primes $p=29,\, 277$ for the
curve with $D=-14.$ Table I below gives the value of $c_p^+(E)$ mod
$p^3$ for both curves and our $p$ in the range $ 5 \leq p < 1000$, while
Table II gives the analogous data for $p$ in the range $11000\,<p\,<
12000.$ Again the values for all our $p$ in the range $p<13500$ can be
found at \cite{tab}.

\smallskip

{\bf Table I}: $c_p^+(E)\cdot p^{-2} \mod p \text{ for }5\leqslant
p<1000\text{ and }p\equiv 1\mod 4 $.

\begin{tabular}{llllll}
$p$& $D=17$&$D=-14$&$p$& $D=17$& $D=-14$\\
& &\\

$ 5 $&$ 3 $&$ 4 $&
 $ 13 $&$ 8 $&$ 4 $\\
 $ 17 $& not valid &$ 7 $&
 $ 29 $&$ 22 $&$ 0 $\\
 $ 37 $&$ 20 $&$ 9 $&
 $ 41 $&$ 29 $&$ 12 $\\
 $ 53 $&$ 45 $&$ 42 $&
 $ 61 $&$ 26 $&$ 60 $\\
 $ 73 $&$ 26 $&$ 56 $&
 $ 89 $&$ 21 $&$ 65 $\\
 $ 97 $&$ 83 $&$ 90 $&
 $ 101 $&$ 59 $&$ 53 $\\
 $ 109 $&$ 34 $&$ 68 $&
 $ 113 $&$ 36 $&$ 47 $\\
 $ 137 $&$ 107 $&$ 126 $&
 $ 149 $&$ 60 $&$ 111 $\\
 $ 157 $&$ 145 $&$ 48 $&
 $ 173 $&$ 44 $&$ 149 $\\
 $ 181 $&$ 70 $&$ 157 $&
 $ 193 $&$ 115 $&$ 11 $\\
 $ 197 $&$ 145 $&$ 54 $&
 $ 229 $&$ 178 $&$ 109 $\\
 $ 233 $&$ 34 $&$ 174 $&
 $ 241 $&$ 141 $&$ 7 $\\
 $ 257 $&$ 199 $&$ 9 $&
 $ 269 $&$ 188 $&$ 139 $\\
 $ 277 $&$ 235 $&$ 0 $&
 $ 281 $&$ 129 $&$ 107 $\\
 $ 293 $&$ 250 $&$ 133 $&
 $ 313 $&$ 69 $&$ 245 $\\
 $ 317 $&$ 237 $&$ 191 $&
 $ 337 $&$ 19 $&$ 151 $\\
 $ 349 $&$ 113 $&$ 263 $&
 $ 353 $&$ 143 $&$ 15 $\\
 $ 373 $&$ 75 $&$ 236 $&
 $ 389 $&$ 257 $&$ 300 $\\
 $ 397 $&$ 78 $&$ 68 $&
 $ 401 $&$ 349 $&$ 340 $\\
 $ 409 $&$ 11 $&$ 313 $&
 $ 421 $&$ 152 $&$ 244 $\\
 $ 433 $&$ 432 $&$ 152 $&
 $ 449 $&$ 423 $&$ 140 $\\
 $ 457 $&$ 288 $&$ 376 $&
 $ 461 $&$ 133 $&$ 37 $\\
 $ 509 $&$ 103 $&$ 407 $&
 $ 521 $&$ 106 $&$ 423 $\\
 $ 541 $&$ 33 $&$ 422 $&
 $ 557 $&$ 276 $&$ 84 $\\
 $ 569 $&$ 423 $&$ 209 $&
 $ 577 $&$ 39 $&$ 212 $\\
 $ 593 $&$ 523 $&$ 18 $&
 $ 601 $&$ 373 $&$ 508 $\\
 $ 613 $&$ 429 $&$ 590 $&
 $ 617 $&$ 133 $&$ 536 $\\
 $ 641 $&$ 285 $&$ 489 $&
 $ 653 $&$ 96 $&$ 540 $\\
 $ 661 $&$ 20 $&$ 330 $&
 $ 673 $&$ 630 $&$ 197 $\\
 $ 677 $&$ 332 $&$ 185 $&
 $ 701 $&$ 105 $&$ 95 $\\
 $ 709 $&$ 437 $&$ 108 $&
 $ 733 $&$ 260 $&$ 462 $\\
 $ 757 $&$ 357 $&$ 672 $&
 $ 761 $&$ 363 $&$ 596 $\\
 $ 769 $&$ 751 $&$ 343 $&
 $ 773 $&$ 13 $&$ 369 $\\
 $ 797 $&$ 123 $&$ 93 $&
 $ 809 $&$ 443 $&$ 212 $\\
\end{tabular}

\begin{tabular}{llllll}
$p$& $D=17$&$D=-14$&$p$& $D=17$& $D=-14$\\
& &\\
 $ 821 $&$ 6 $&$ 347 $&
 $ 829 $&$ 645 $&$ 823 $\\
 $ 853 $&$ 48 $&$ 635 $&
 $ 857 $&$ 5 $&$ 502 $\\
 $ 877 $&$ 132 $&$ 603 $&
 $ 881 $&$ 82 $&$ 591 $\\
 $ 929 $&$ 845 $&$ 766 $&
 $ 937 $&$ 341 $&$ 100 $\\
 $ 941 $&$ 253 $&$ 642 $&
 $ 953 $&$ 794 $&$ 866 $\\
 $ 977 $&$ 548 $&$ 98 $&
 $ 997 $&$ 302 $&$ 401 $\\
\end{tabular}

\smallskip

{\bf Table II}: $c_p^+(E)\cdot p^{-2} \mod p\text{ for
}11000<p<12000\text{ and }p\equiv 1\mod 4 $.

\begin{tabular}{llllll}
$p$& $D=17$&$D=-14$&$p$& $D=17$& $D=-14$\\
& & & & &\\
$ 11057 $&$ 3236 $&$ 10336 $&
 $ 11069 $&$ 7768 $&$ 6637 $\\
 $ 11093 $&$ 9234 $&$ 5437 $&
 $ 11113 $&$ 832 $&$ 9242 $\\
 $ 11117 $&$ 6204 $&$ 7965 $&
 $ 11149 $&$ 8885 $&$ 1364 $\\
 $ 11161 $&$ 1292 $&$ 1636 $&
 $ 11173 $&$ 587 $&$ 10503 $\\
 $ 11177 $&$ 6184 $&$ 4427 $&
 $ 11197 $&$ 8804 $&$ 6750 $\\
 $ 11213 $&$ 6409 $&$ 8508 $&
 $ 11257 $&$ 192 $&$ 1839 $\\
 $ 11261 $&$ 700 $&$ 6850 $&
 $ 11273 $&$ 5932 $&$ 510 $\\
 $ 11317 $&$ 1969 $&$ 2892 $&
 $ 11321 $&$ 5451 $&$ 10402 $\\
 $ 11329 $&$ 5635 $&$ 9145 $&
 $ 11353 $&$ 3322 $&$ 7820 $\\
 $ 11369 $&$ 6790 $&$ 11276 $&
 $ 11393 $&$ 4532 $&$ 358 $\\
 $ 11437 $&$ 10570 $&$ 3120 $&
 $ 11489 $&$ 8715 $&$ 10941 $\\
 $ 11497 $&$ 4837 $&$ 6424 $&
 $ 11549 $&$ 7265 $&$ 2757 $\\
 $ 11593 $&$ 225 $&$ 369 $&
 $ 11597 $&$ 8864 $&$ 7113 $\\
 $ 11617 $&$ 10691 $&$ 1052 $&
 $ 11621 $&$ 7500 $&$ 6521 $\\
 $ 11633 $&$ 293 $&$ 5463 $&
 $ 11657 $&$ 10665 $&$ 4770 $\\
 $ 11677 $&$ 10365 $&$ 11566 $&
 $ 11681 $&$ 6023 $&$ 5351 $\\
 $ 11689 $&$ 11553 $&$ 3152 $&
 $ 11701 $&$ 5851 $&$ 11618 $\\
 $ 11717 $&$ 10185 $&$ 8521 $&
 $ 11777 $&$ 10882 $&$ 3487 $\\
 $ 11789 $&$ 6221 $&$ 3509 $&
 $ 11801 $&$ 10632 $&$ 3148 $\\
 $ 11813 $&$ 2123 $&$ 3767 $&
 $ 11821 $&$ 7340 $&$ 128 $\\
 $ 11833 $&$ 1715 $&$ 9412 $&
 $ 11897 $&$ 8766 $&$ 10281 $\\
 $ 11909 $&$ 6032 $&$ 11519 $&
 $ 11933 $&$ 1190 $&$ 1783 $\\
 $ 11941 $&$ 5023 $&$ 6379 $&
 $ 11953 $&$ 10988 $&$ 1162 $\\
 $ 11969 $&$ 11669 $&$ 11573 $&
 $ 11981 $&$ 1742 $&$ 8384 $\\
\end{tabular}

\medskip

\n Finally, for the curve $y^2=x^3+14x$ and the two exceptional primes
$p=29,277$, we have
\begin{align*}
c_{29}^+(E)&\equiv 27\cdot
29^3\mod 29^4,\\
c_{277}^+(E)&\equiv 155\cdot 277^3\mod 277^4.
\end{align*}

\n John Coates\hb
\n DPMMS, University of Cambridge\hb
\n Centre for Mathematical Sciences\hb
\n Wilberforce Road\hb
\n Cambridge CB3 0WB, England.\hb
\n e-mail: J.H.Coates@dpmms.cam.ac.uk

\bigskip

\n Zhibin Liang\hb
\n School of Mathematical Sciences\hb
\n Capital Normal University\hb
\n Beijing, China.\hb
\n e-mail: liangzhb@gmail.com

\bigskip

\n R. Sujatha\hb
\n School of Mathematics\hb
\n Tata Institute of Fundamental Research\hb
\n Homi Bhabha Road\hb
\n Bombay 400 005, India.\hb
\n e-mail: sujatha@math.tifr.res.in


\begin{thebibliography}{999}

\bibitem{BGS} {\sc D. Bernardi, C. Goldstein, N. Stephens},
\newblock{\it Notes $p$-adiques sur les
courbes elliptiques}, \newblock{Crelle Jour.}  {\bf 351} (1984),
129-170.

\bibitem{Ca} {\sc J. Cassels}, \newblock{\it A note on the division
values of $p(u)$}, \newblock {Proc. Camb. Phil. Soc.} {\bf 45} (1949),
167-172.


\bibitem{C}{\sc J. Coates}, \newblock{\it Elliptic curves with complex
multiplication and Iwasawa theory}, \newblock {Bull. London
Math. Soc.} {\bf 23} (1991), 321-350.

\bibitem{CS}{\sc J. Coates, R. Sujatha}, \newblock{Elliptic curves
with complex multiplication and $L$-values}, Book in preparation.

\bibitem{CW}{\sc J. Coates, A. Wiles}, \newblock{\it On the conjecture
of Birch and Swinnerton-Dyer}, \newblock {Invent. Math.} {\bf 39} (1977
), 223-251.

\bibitem{CW1}{\sc J. Coates, A. Wiles}, \newblock{\it On $p$-adic
$L$-functions and elliptic units}, \newblock {Jour. Australian
Math.Soc.} {\bf 26} (1978), 1-25.

\bibitem{DD}{\sc T. Dokchitser, V. Dokchitser,} {\it On the Birch
Swinnerton-Dyer quotients modulo squares}, \newblock{Ann. Math.} (To
appear).

\bibitem{Gr}{\sc B. Gross}, \newblock{Arithmetic on elliptic
curves with complex multiplication}, \newblock Lecture Notes in
Math. {\bf 776}, Springer (1980).

\bibitem{GS}{\sc C. Goldstein, N. Schappacher}, \newblock{\it S\'eries
d'Eisenstein et fonctions $L$ de courbes elliptiques \`a
multiplication complexe}, \newblock {Crelle J.} {\bf 327} (1981),
184-218.

\bibitem{FK}{\sc T. Fukuda, K. Komatsu} \newblock {\it ${\Bbb 
Z}_p$-extensions 
associated to elliptic curves
with complex multiplication}, \newblock {Math. Proc. Camb. Phil. Soc.} 
{\bf 137} (2004), 541-550.

\bibitem{Ko} {\sc V. Kolyvagin}, \newblock{\it Euler systems}, in
\newblock{The
Grothendieck
Festschrift}, Vol. II,  435--483, Progr. Math. {\bf  87}, Birkhäuser
Boston, Boston, MA, 1990.

\bibitem{Ne}{\sc J. Nekov\'a\v r,}\newblock{\em Selmer
complexes}, Ast\'erisque {\bf 310} (2007).

\bibitem{PR}{\sc B. Perrin-Riou}, \newblock{\it Arithm\'etique des
courbes elliptiques et th\'eorie d'Iwasawa}, \newblock
{Bull. Soc. Math. France Suppl. M\'emoire} {\bf 17} (1984).

\bibitem{BPR}{\sc D. Bernardi, B. Perrin-Riou}, \newblock{\it
Variante $p$-adique de la conjecture de Birch et Swinnerton-Dyer (le cas
supersingulier)}, \newblock{C. R. Acad. Sci. Paris Sér. I
Math.} {\bf 317}
(1993), 227-232.

\bibitem{PR1}{\sc B. Perrin-Riou}, \newblock{\it
Arithm\'e©tique des
courbes
elliptiques \`a r\'eduction super singuli\'ere en $p$},
\newblock{Experiment.
Math.} {\bf  12}(2003), 155-186.

\bibitem{Ru}{\sc K. Rubin}, \newblock{\it The ``main conjectures'' of
Iwasawa theory for imaginary quadratic fields}, \newblock
{Invent. Math.} {\bf 103} (1991), 25-68.


\bibitem{SW} {\sc W. Stein, C. Wuthrich}, \newblock{\it Computations
about Tate-Shafarevich groups using Iwasawa theory}, Preprint (2008).

\bibitem{Se} {\sc J.-P. Serre}, \newblock {Corps Locaux}, \newblock GTM
{\bf 67}, Springer (1979).

\bibitem{tab} URL: http://www.cnu.edu.cn/mathpage/upload/liangzibin/thesis/20081104074048.rar

\bibitem{vdw} {\sc B. Van der Waerden}, \newblock{Modern Algebra},
Vols. I and II, Ungar (1953).

\bibitem{We} {\sc A. Weil}, \newblock{Elliptic functions according to
Eisenstein and Kronecker}, \newblock Springer (1976).

\bibitem{AT} {\sc J. Tate}, \newblock{On the conjecture of Birch and
 Swinnerton-Dyer and a geometric analog}, \newblock{S\'eminaire Bourbaki}, {\bf
 306} (1966).


\end{thebibliography}
\end{document}